\magnification=\magstep1
\input amssym.def
\input amssym
     
\baselineskip=12pt
\hsize=6.3truein
\vsize=8.7truein
\font\footsc=cmcsc10 at 8truept
\font\footbf=cmbx10 at 8truept
\font\footrm=cmr10 at 10truept
\footline={\footsc the electronic journal of combinatorics
        {\footbf 12} (2005)  \#R21\hfil\footrm\folio}
\font\twelverm=cmr12
%
\def\parno{\par\noindent}
\def\sr{\scriptscriptstyle}
\def\vb{\vrule height 14pt depth 7pt} 
\def\ts{\tabskip0pt}	
\def\ss1{\noalign{\vskip -1pt}}
\def\vs1{\noalign{\vskip-5pt}}
\def\gb{{\goth b}}
\def\gg{{\goth g}}
\def\gh{{\goth h}}
\def\gl{{\goth l}}
\def\gm{{\goth m}}
\def\gn{{\goth n}}
\def\gs{{\goth s}}
\def\bBo{{\bf B}}
\def\bG{{\bf G}}
\def\bH{{\bf H}}
\def\bS{{\bf S}}
\def\bD{{\bf D}}
\def\bT{{\bf T}}
\def\a{\alpha}

\def\l{\lambda}
\def\om{\omega}
\def\Na{{\Bbb N}}  
\def\Dscr{{\cal D}}
\def\Oscr{{\cal O}}
\def\Cscr{{\cal C}}
\def\Bscr{{\cal B}}
\def\Vscr{{\cal V}}
\def\sh{{\rm sh}\,}
\def\Lie{{\rm Lie}\,}
\def\qed{\hbox{\hskip 4pt
                \vrule width 5pt height 6pt depth 1.5pt\hskip 2pt}}
\def\QED{\par\hfill\qed \par}
\def\pr{^{\prime}}
\def\prpr{^{\prime\prime}}
\def\st{\subset}
\def\steq{\subseteq}
\def\ra{\rightarrow}
\def\Ra{\Rightarrow}
\def\ov{\overline}

\def\la{\leftarrow}
\def\pmat{\pmatrix}
\font\mas=msbm10
\def\bR{\hbox{\mas\char'124}}
\def\Pf{\goodbreak\noindent{\bf Proof.}\parno}
\def\do{\buildrel \sr D \over \leq}     
\def\dg{\buildrel \sr D \over \geq}
\def\dos{\buildrel \sr D \over <}
\def\dgs{\buildrel \sr D \over >}
\def\go{\buildrel \sr G \over \leq}     

\def\cho{\buildrel \sr C \over \leq}     
\def\chos{\buildrel \sr C \over <}
\def\chgs{\buildrel \sr C \over >}
\def\chg{\buildrel \sr C \over \geq}
\def\La{\Leftarrow}
%
\def\<#1>{\langle#1\rangle}
\overfullrule=0pt
%
\font\bigrm=cmr12 at 14pt
\font\bigbf=cmb12 at 14pt

\centerline{\bigbf The combinatorics of orbital varieties closures}
\vskip 1.3mm
\centerline{\bigbf of nilpotent order 2 in sl${}_n$}

\bigskip\bigskip

\centerline {{\bigrm Anna Melnikov}\footnote {*}{Supported in part by the Minerva Foundation, Germany,
Grant No. 8466}}
\smallskip
\centerline {Department of Mathematics,}
\centerline {University of Haifa,}
\centerline {31905 Haifa, Israel}
\centerline {and}
\centerline {Department of Mathematics,}
\centerline {the Weizmann Institute of Science,}
\centerline {76100 Rehovot, Israel}
\centerline {\tt melnikov@math.haifa.ac.il	}

\bigskip

\centerline{\footrm Submitted: Sep 12, 2002; Accepted: Apr 28, 2005; Published: May 6, 2005}
\centerline{\footrm Mathematics Subject Classifications: 05E10, 17B10}

\bigskip\bigskip

\centerline {\bf Abstract.}
\smallskip
{\narrower\noindent We consider two partial orders
on the set of standard Young tableaux. The first one is induced
to this set  from the weak right
order on symmetric group by Robinson-Schensted algorithm.
The second one is induced to it from the dominance order
on Young diagrams by considering
a Young tableau as a chain of Young diagrams. We prove
that these two orders of completely different nature coincide
on the subset of Young tableaux with 2 columns or with 2 rows.  
This fact has very interesting geometric implications for orbital 
varieties of nilpotent order 2 in special linear algebra $sl_n.$\par}

\bigskip

{\bf 1. \ \ Introduction}
\parno
{\bf 1.1}\ \ \ Let $\bS_n$ be a symmetric group, that is a group of permutations of  $\{1,2,\ldots,n\}$.
Respectively, let ${\Bbb S}_n$ be a group of permutations
of $n$ positive integers $\{m_1<m_2<\ldots< m_n\}$ 
where $m_i\geq i.$ It is obvious that there is a bijection from ${\Bbb S}_n$
onto $\bS_n$ obtained by $m_i\ra i,$ so we will use the notation ${\Bbb S}_n$ in all the cases
where the results apply to both ${\Bbb S}_n$ and  $\bS_n.$

In this paper we write a permutation in a word form
$$ w =[a_{\sr 1}, a_{\sr 2},\ldots, a_n]\ ,\quad {\rm where}
                                    \ a_i=w(m_i).\eqno{(*)}.$$
All the words considered in this paper are permutations, i.e. with distinct letters only.

Set $p_w(m_i):=j$ if $a_j=m_i,$ in other words,
$p_w(m_i)$ is the place (index) of $m_i$ in the word form
of $w.$ (If $w\in \bS_n$ then $p_w(i)=w^{-1}(i).$)
\par
We consider the right weak (Bruhat) order on ${\Bbb S}_n$ that is 
we put $w\do y$ if for all 
$i,j\ :\ 1\leq i<j\leq n$ the condition $p_w(m_j)<p_w(m_i)$ implies
$p_y(m_j)<p_y(m_i).$ Note that $[m_{\sr 1},m_{\sr 2},\ldots,m_n]$ is the minimal word 
and $[m_n,m_{n-\sr 1},\ldots,m_{\sr 1}]$ is the maximal word in this 
order. 
\parno
{\bf 1.2}\ \ \ 
Let\ $\l=(\l_1\geq \l_2 \geq \cdots \geq \l_k > 0)$
be a partition of $n$ and
$\l\pr:=(\l\pr_1\geq \l\pr_2\geq \cdots\geq \l\pr_l>0)$
the conjugate partition, that is $\l\pr_i=
\sharp\{j\ |\ \l_j\geq i\}.$ In particular, 
$\lambda^{\prime}_1=k.$
\par
We define the corresponding Young
diagram $D_\l$  of $\l$  to be an array of $k$  columns of boxes
starting from the top  with the $i$-th column containing 
$\l_i$ boxes. Note that it is more customary that $\l$ defines the rows of the
diagram and $\l\pr$ defines the columns, but in the present context
we prefer this convention for the simplicity of notation.
Let $\bD_n$ denote the set of all Young diagrams with $n$ boxes.
\par
We use the dominance order on partitions. It is a partial order defined
as follows. 
Let  $\l=(\l_{\sr 1},\cdots,\l_k)$ and 
$\mu = (\mu_{\sr 1},\cdots,\mu_j)$ be partitions of $n.$  
Set
$\l \geq \mu$ if for each
$i\ :\ 1\leq i\leq \min(j,k)$ one has
$$ \sum^i_{m=1}\l_m \geq \sum^i_{m=1}\mu_m \ . $$
\parno
{\bf 1.3}\ \ \ 
Fill the boxes of the Young diagram
$D_\l$ with $n$ distinct positive
integers $m_1<m_2<\ldots<m_n.$  If the entries increase in rows from left to right and in
columns from top to bottom, we call such an array a Young tableau or
simply a tableau.
If the numbers in a
tableau form the set of integers from 1 to $n$,
the tableau is called standard. 
\par
Let ${\Bbb T}_n$ denote the set of tableaux with 
$n$ positive entries $\{m_1<m_2<\ldots< m_n\}$ 
where $m_i\geq i,$
and respectively let $\bT_n$ denote the set of standard tableaux.
Again, the bijection from ${\Bbb T}_n$
onto $\bT_n$ is obtained by $m_i\ra i,$ and we will use the notation ${\Bbb T}_n$ in all the cases
where the results apply to both ${\Bbb T}_n$ and  $\bT_n.$  
The Robinson-Schensted algorithm (cf. [Sa,\S 3], or  [Kn, 5.1.4], or [F, 4.1] ) gives the bijection $w\mapsto (T(w),Q(w))$ from ${\Bbb S}_n$ onto the set of pairs of tableaux of the same 
shape.  For each $T\in {\Bbb T}_n$  set 
$\Cscr_T=\{w\ |\ T(w)=T\}.$ It is called a Young cell. The
right weak order on ${\Bbb S}_n$ induces a natural order relation $\do$
on ${\Bbb T}_n$ as follows. We say that
$T\do S$ if there exists a sequence of tableaux $T=P_1,\ldots, P_k=S$ 
such that for each
$j\ :\ 1\leq j<k$ there exists a pair $w\in \Cscr_{P_j},\ y\in 
\Cscr_{P_{j+1}}$ satisfying $w\do y.$ 
\par
I would like to explain the notation $\do.$ I use it in honor of M. Duflo 
who was the first to discover  the implication of the weak  
order on Weyl group for the primitive spectrum of the corresponding enveloping algebra
(cf. [D]). I would like to use the
notation since his result was the source of my personal interest to the different combinatorial 
orderings of Young tableaux.
\par
Consider $\bS_n$ as a Weyl group of $\gs\gl_n(\Bbb C).$ By Duflo, there is a surjection from $\bS_n$
onto the set of primitive ideals (with infinitesimal character).   Let us define
the corresponding primitive ideal by $I_w.$ By [D], $w\do y$ implies $I_w\subseteq I_y.$
As it was shown by A. Joseph [J], $I_w$ and $I_y$ coincide
iff $w$ and $y$ are in the same Young cell. Together these two facts show that 
the order $\do$ is well defined on $\bT_n.$
\par
As shown in [M1, 4.3.1], 
one may have
$T,S\in \bT_n$ for which $T\dos S;$ yet for any
$w\in \Cscr(T),\ y\in \Cscr(S)$ one has $w\not\dos y.$ Thus, it is essential to define it through
the sequence of tableaux.
\parno
{\bf 1.4}\ \ \ Take $T\in {\Bbb T}_n$ and let 
$\sh(T)$ be the underlying diagram of $T.$ We will write it
as $\sh(T)=(\l_1,\ldots,\l_k)$ where $\l_i$ is the length
of the $i-$th column. 
Given  $i,j\ : 1\leq i<j\leq n$ we define 
$\pi_{i,j}(T)$ to be the tableau obtained from $T$ by removing $m_{\sr 1},\ldots,m_{i-\sr 1}$
and $m_{j+\sr 1},\ldots,m_n$ by ``jeu de taquin'' (cf. [Sch] or 2.10).
Put $D_{\<i,j>}(T):=\sh (\pi_{i,j}(T)).$ We define the following
partial order on ${\Bbb T}_n$ which we call the chain order.
We set $T\cho S$ if for any $i,j\ :\ 1\leq i<j\leq n$
one has $D_{\<i,j>}(T)\leq D_{\<i,j>}(S).$
\par
This order is obviously well defined.
\parno
{\bf 1.5}\ \ \ The above constructions give two purely 
combinatorial orders on ${\Bbb T}_n$
which are moreover of an entirely different nature. 

Given two partial orders $\buildrel  a \over \leq$ and $\buildrel  b \over \leq$
on the same set S, call $\buildrel  b \over \leq$ an extension of 
$\buildrel  a \over \leq$ if $s{ \buildrel  a \over \leq}t$ implies
$s{ \buildrel  b \over \leq}t$  for any $s,t\in S.$

As we explain in 1.11,
$\cho$ is an extension of $\do$ on $\bT_n.$ Moreover, 
 these two orders coincide for $n\leq 5$ and 
$\cho$ is a proper extension of $\do$ for $n\geq 6,$ as shown in [M].
\par
There is a significant simplification when one considers only
tableaux with two columns. Let us denote the subset
of tableaux with two columns in ${\Bbb T}_n$ by ${\Bbb T}^2_n.$
We show that for $S,T\in {\Bbb T}_n^2$ one has 
$T\cho S$ if and only if $T\do S.$ 
Moreover, for any $T\in{\Bbb T}_n^2$ we construct a canonical
representative $w_{\sr T}\in \Cscr_T$ such that
$T\chos S$ if and only if $w_{\sr T}\dos w_{\sr S}.$
\parno
{\bf 1.6}\ \ \ Given a set $S$ and a partial order $\buildrel  a \over \leq,$
the cover of $t\in S$  in this order is the set of all $s\in S$ such
that
$t{\buildrel a\over <} s$ and there is no $p\in S$ such that $t{\buildrel  a \over <}p{\buildrel  a \over 
<} s.$ We will denote it by $\Dscr_a(t).$
 
As explained in [M1], in general, even an 
inductive description of  $\Dscr_D(T)$ is a very complex task. 
Yet, in 3.16 we provide the exact description of $\Dscr_D(T)$ (which is a cover
in $\cho$ as well) for any $T\in {\Bbb T}_n^2.$
\parno
{\bf 1.7}\ \ \ For each tableau $T$ let $T^\dagger$ denote
the transposed
tableau. Obviously, $T\chos S$
iff $S^\dagger\chos T^\dagger.$ By Schensted-Sch\"utzenberger
theorem (cf. 2.14), it is obvious that
 $T\dos S$
iff $S^\dagger\dos T^\dagger.$  Consequently, the above results
can be translated to tableaux with two rows. 
\parno
{\bf 1.8}\ \ \ Let us finish the introduction by explaining why these two 
orders are of interest and what implication our results have
for the theory of orbital varieties.
\par
Orbital varieties arose from
the works of N. Spaltenstein ([Sp1] and [Sp2]), 
and R. Steinberg ([St1] and [St2]) during their studies of the 
unipotent 
variety of a semisimple group $\bG.$ 
\par
Orbital varieties are 
the translation of these components from the unipotent
variety of $\bG$ to the nilpotent cone of $\gg =\Lie(\bG).$
They are defined as follows.
\par
Let $\bG$ be a connected semisimple 
finite dimensional
complex algebraic group. Let $\gg$ be its Lie algebra  and
$U(\gg)$ be the enveloping algebra of $\gg.$
Consider the adjoint action of $\bG$ on $\gg.$
Fix some triangular decomposition  $\gg=\gn\bigoplus
\gh\bigoplus\gn^-.$  A $\bG$ orbit $\Oscr$ in $\gg$ is called
nilpotent if it consists of nilpotent elements, that is if
$\Oscr=\bG\, x$ for some $x\in\gn.$ 
The intersection $\Oscr\cap\gn$ is reducible. Its irreducible components
are called orbital varieties associated to $\Oscr.$
They are Lagrangian subvarieties of $\Oscr.$ According to the orbit method philosophy, they should play an important role in the representation theory of corresponding Lie algebras. Indeed, they play the key role 
in the study of
primitive ideals in $U(\gg).$ They also play an important
role in Springer's Weyl group representations
described in terms of fixed point sets $\Bscr_u$ where $u$ is a unipotent element acting on
the flag variety $\Bscr.$
\par
Orbital varieties are very interesting objects
from the point of view of algebraic geometry. 
Given an orbital variety $\Vscr$, one can easily find 
the nilradical $\gm_{\sr \Vscr}$ of a standard parabolic 
subalgebra of the smallest dimension 
containing $\Vscr$. Consider an orbital variety 
closure as an algebraic variety in the affine linear space  $\gm_{\sr \Vscr}.$
Then the vast majority of orbital varieties are not complete 
intersections.  
So, orbital varieties are examples of algebraic varieties
which are both Lagrangian subvarieties and not 
complete intersections.
\parno
{\bf 1.9}\ \ \ There are many hard open questions involving orbital varieties.
Their only general description
was given by R. Steinberg [St1]. Let us explain it briefly.
\par
Let $R\st\gh^*$ denote the set of 
non-zero roots,  $R^+$ the set of positive roots corresponding 
to $\gn$ and $\Pi\st R^+$ the resulting set of 
simple roots. Let $W$ be the Weyl group for the pair $(\gg,\ \gh).$ 
For any $\a \in R$ let $X_\a$ be the corresponding root space.
\par
For $S,S\pr\st R$ and $w\in W$ set $S\cap^w S\pr:=\{\a\in S\ :\ 
\a\in w(S\pr)\}.$ Then set
$$\gn\cap^w\gn:=\bigoplus\limits_{\a\in R^+\cap^w R^+}X_{\a}.$$
This is a subspace of $\gn.$ For each closed irreducible  subgroup  $\bH$ of $\bG$ let
$\bH(\gn\cap^w\gn)$ be the set of $\bH$
conjugates of $\gn\cap^w\gn.$ It is an irreducible locally 
closed subvariety. Let $\ov *$ denote the (Zariski) closure
of a variety $*.$
\par
Since there are only finitely many nilpotent 
orbits in $\gg$, it follows that there exists a unique 
nilpotent orbit which we denote by $\Oscr_w$ 
such that 
$\ov{\bG(\gn\cap^w\gn)}=\ov\Oscr_w.$
\par
Let $\bBo$ be the standard
Borel subgroup of $\bG$, i.e.  such that 
$\Lie(\bBo)=\gb=\gh\bigoplus\gn.$
A result of Steinberg [St1] asserts that 
$\Vscr_w:=\ov{\bBo(\gn\cap^w\gn)}\cap\Oscr_w$ is an orbital
variety and that the map
$\varphi:w\mapsto \Vscr_w$
is a surjection of $W$ onto the set of orbital varieties.
The fibers of this mapping, namely $\varphi^{-1}(\Vscr)=
\{w\in W\ :\ 
\Vscr_w=\Vscr\}$ are called geometric cells. 
\par
This description is not very satisfactory 
from the geometric 
point of view since a $\bBo$ invariant subvariety
generated by a linear space is a very complex object.
For example, one can describe the regular functions
(differential operators) on $\ov\Vscr_w$ or on $\Vscr_w$
only in some special cases. 
\parno
{\bf 1.10}\ \ \ On the other hand, there exists a very nice combinatorial
characterization of orbital varieties in  $\gs\gl_n$ in terms of
Young tableaux. Indeed, in that case $\Vscr_w$ and $\Vscr_y$ coincide
iff $w$ and $y$ are in the same Young cell. Moreover, let $\Oscr_w=\bG\Vscr_w$ be the 
corresponding nilpotent orbit, then its Jordan form is defined by
$\mu=(\sh T_w)\pr.$ Let us denote such orbit by $\Oscr_{\mu}.$
\par
Recall the order relation on Young diagrams from 1.2.
A result of Gerstenhaber (see [H, \S 3.10] for example)
describes the closure of a nilpotent orbit.
\proclaim Theorem. Let $\mu$ be a partition of $n.$ One has
$$\ov\Oscr_{\mu}=\coprod_{\l\vert \l\geq \mu}\Oscr_{\l}$$
\par
\parno
{\bf 1.11}\ \ \ 
Define geometric order on $\bT_n$ by $T\go S$ if 
$\ov\Vscr_S\st \ov\Vscr_T.$ In general, the combinatorial 
description of this order is an open (and very difficult) task.
On the other hand, both $\do$ and $\cho$ are connected to
$\go$ as follows.
\par
Let us identify $\gn$ with the subalgebra of strictly
upper-triangular matrices.  Any $\a\in R^+$ can be decomposed into the sum of 
simple roots $\a=\sum_{k=i}^{j-1}\a_k$ where $i< j.$ 
Then the root space $X_\a$
is identified with $X_{i,j}.$
By [JM, 2.3], $X_{i,j}\in\gn\cap^w\gn$ if and only if 
$p_w(i)<p_w(j).$
Thus, $w\do y$ implies $\gn\cap^y\gn\st \gn\cap^w\gn$, hence, also 
$\ov\Vscr_y\subset \ov\Vscr_w$ and $\ov\Oscr_y\subset\ov\Oscr_w.$
Therefore, $\go$ is an extension of $\do$ on $\bT_n.$  
\par
On the other hand, note that $T\go S$ implies, in particular, the inclusion of 
corresponding orbit closures so that (via Gerstenhaber's construction)
$T\go S$ implies $\sh(T)\leq \sh(S).$ As shown in [M1, 4.1.1], the projections
on the Levi factor of standard parabolic subalgebras of $\gg$ 
preserve orbital variety closures. Moreover, in the case
of $\gs\gl_n$ one has $\pi_{i,j}(\ov\Vscr_T)=
\ov\Vscr_{\pi_{i,j}(T)}$ for any $i,j\ :\ 1\leq i<j\leq n$ where $\pi_{i,j}(T)$ is obtained from $T$
by jeu de taquin and $\Vscr_{\pi_{i,j}(T)}$
is an orbital variety in the corresponding Levi factor. 
Thus,  $T\go S$
implies $\pi_{i,j}(T)\go \pi_{i,j}(S).$ Altogether, this provides
that $\cho$ is an extension of $\go.$
\par
Consequently, $\cho$ is an extension of $\go$ and 
$\go$ is an extension of $\do.$ 
All three orders coincide for $n\leq 5,$ and 
$\cho$ is a proper extension of
$\go$ which is, in turn, a proper extension of $\do$ for $n\geq 6$ as shown in [M].
\par
However, our results  show that $\do$ and $\cho$ coincide on $\bT^2_n$ and there 
they provide a full combinatorial
description of $\go.$ 
\par
Consider $\Vscr_T$ where $T\in\bT_n^2.$ 
For any $X\in\Vscr_T$
one has $X\in\Oscr_{\sh(T)}$, that is $X$ is an element of
nilpotent order $2$ or in other words $X^2=0.$ Thus, we get a complete combinatorial
description of inclusion of orbital varieties closures of nilpotent order $2$ in $\gs\gl_n.$ 
\parno
{\bf 1.12}\ \ \ The body of the paper consists of two sections.
\par
In section 2 we explain all the background in combinatorics of 
Young tableaux
essential in  the subsequent analysis and set the notation. 
In particular, we explain Robinson-Schensted insertion from the left
and jeu de taquin. I hope this part makes the paper 
self-contained.
\par
In section 3  we work out the machinery for 
comparing $\do$ and $\cho$ and show that they coincide.
The main technical result of the paper is stated in 3.5 and proved in 3.11.
Further in 3.12, 3.13 and 3.14 we explain the implications of this result
for $\do,\ \go$ and $\cho.$ In 3.16 we give the exact description of
$\Dscr_{G}(T)$ for $T\in \bT_n^2.$ Finally, in 3.17 we explain the corresponding
facts for the tableaux with two rows.
\par

\bigskip
{\bf 2. \ \ Combinatorics of Young tableaux}
\parno
{\bf 2.1}\ \ \ 
Recall from 1.1 $(*)$ the presentation of $w\in\bS_n$  in 
the word form. Given $w\in \bS_n$, set 
$$\tau(w):=\{i\ :\ p_w(i+1)<p_w(i)\},$$
that is $\tau(w)$ is the set of left descents of $w.$

Note that if $w\do y$ then $\tau(w)\steq \tau(y).$
\parno
{\bf 2.2}\ \ \ Given a word or a tableau  $*$, we denote by $\<*>$
the set of its entries. Introduce the following useful notational conventions.
\item {(i)} For $m\in \<w>$ set $w\setminus\{m\}$ to be the word obtained from $w$ by deleting
$m,$ that is if $m=a_i$ then $w\setminus \{m\}:= [a_{\sr 1},\ldots,a_{i-\sr 1},a_{i+\sr 1},\ldots,a_n].$
\item {(ii)} For the words $x=[a_{\sr 1},\ldots,a_n]$ and $y=[b_{\sr 1},\ldots,b_m]$ such that
 $\<x>\cap \<y>=\emptyset$ we define a colligation $[x,y]:=[a_{\sr 1},\ldots,a_n,b_{\sr 1},\ldots,b_m]$.
\item {(iii)} For a word $w=[a_{\sr 1},\ldots,a_n]$ set 
$\ov w$ to be the word with reverse order, that is 
$\ov w:=[a_n,a_{n-\sr 1},\ldots,a_{\sr 1}].$
\par
Given $i,j\ :\ 1\leq i<j\leq n$, set ${\Bbb S}_{\<i,j>}$
to be a (symmetric) group of permutations of $\{m_k\}_{k=i}^j$.
Let us define 
projection  $\pi_{i,j}:{\Bbb S}_n\ra {\Bbb S}_{\<i,j>}$ by 
omitting  all the letters $m_1,\ldots,m_{i-1}$
and $m_{j+1},\ldots,m_n$ 
from word $w\in {\Bbb S}_n$, 
i.e. 
$\pi_{i,j}(w)=w\setminus\{m_1,\ldots,m_{i-1},m_{j+1},\ldots,m_n\}.$ For $w\in\bS_n$ it  
is obvious that $\tau(\pi_{i,j}(w))=\tau(w)\cap\{k\}_{k=i}^{j-1}.$

\proclaim Lemma. Let $w,y$ be in ${\Bbb S}_n.$
\item{(i)} 
For any $a\not \in\{m_i\}_{i=1}^n$
one has $w\do y$ iff $[a,w]\do [a,y].$
\item {(ii)} For $w,y$ such that 
$\pi_{1,n-1}(y)=\pi_{1,n-1}(w)$
and $p_w(m_n)=1,\ p_y(m_n)>1$ one has $w\dgs y.$
\item {(iii)} $w\dos y$ iff $\ov y\dos \ov w.$
\item{(iv)} If $w\do y$ then
$\pi_{i,j}(w)\do \pi_{i,j}(y)$ for any $i,j\ : 1\leq i<j\leq n.$ 
\par
All four parts of the lemma are obvious.
\parno
{\bf 2.3}\ \ \ We will use the following notation for
tableaux. 
Let $T$ be a tableau and let $T^i_j$ for $i,j\in \Na$ denote the entry on the intersection of
the $i$-th row and the $j$-th column.
Given $u$ an entry of $T$, set $r_{\sr T}(u)$ to be the
number of the row, $u$ belongs to and $c_{\sr T}(u)$ to
be the number of the column, $u$ belongs to. Set 
$$\tau(T):=\{i\ :\ r_{\sr T}(i+1)>r_{\sr T}(i)\}.$$
Let $T_i$ denote the $i$-th column of $T.$ Let  $\om_i(T)$   
denote the largest entry of $T_i$. 
\par
We  consider a tableau as a  matrix $T:=(T_i^j)$  
and  write $T$ by columns: $T = (T_1, \cdots, T_l)$
\par
For $i,j\ :\ 1\leq i<j\leq l$ we set $T_{i,j}$  to be a 
subtableau of $T$ consisting of columns from $i$ 
to $j$, that is $ T_{i,j}=(T_i,\cdots, T_j).$
For each tableau $T$ let $T^{\dagger}$ denote the transposed tableau.
Note that $\sh(T^\dagger)=\sh(T)\pr.$
\parno
{\bf 2.4}\ \ \ Given $D_\l \in \bD_n$ with $\l=(\l_1,\cdots,\l_j)$,
we define a corner box (or simply, a corner) of the Young diagram 
to be a box with
no neighbours to right and below.
\par
For example, in $D$ below all the corner boxes are labeled by $X$.
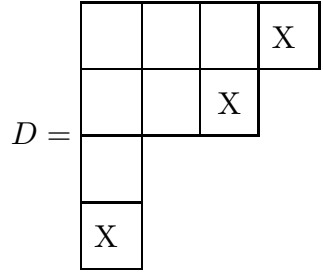
$$ D =
\vcenter{
\halign{& \hfill#\hfill
\tabskip4pt\cr
\multispan{9}{\hrulefill}\cr
\ss1
\vb & \quad & \vb &  \quad & \vb  & \quad & \vb & X &\ts\vb\cr
\vs1
\multispan{9}{\hrulefill}\cr
\ss1
\vb & \ & \vb &  \ & \vb  &X&\ts\vb\cr
\vs1
\multispan{7}{\hrulefill}\cr
\ss1
\vb & \ & \ts\vb \cr
\vs1
\multispan{3}{\hrulefill}\cr
\ss1
\vb  &X &\ts\vb\cr
\vs1
\multispan{3}{\hrulefill}\cr}} $$
The entry of a tableau in a corner is called a corner entry.
Take  $D_\l$ with  $\l=(\l_1,\cdots,\l_k).$ Then  there is a 
corner  
entry $\om_i(T)$ at the corner  $c$ with coordinates $(\l_i,i)$ iff 
$\l_{i+1}<\l_i.$
\parno
{\bf 2.5}\ \ \ We now define the insertion algorithm. 
Consider a column  $C = \pmat{a_1\cr \vdots\cr a_r\cr}.$ 
Given $j \in {\Na}^+\setminus \<C>$, 
let $a_i$  be the  smallest entry
greater then $j,$ if exists.  
Set
$$j\ra C:= \cases{
\left( \matrix{a_{\sr_1}\cr\vdots\cr a_{i-\sr 1}\cr j\cr a_{i+\sr 1}\cr\vdots\cr} \right),\ \ 
j_{\sr C}=a_i & if $j<a_r$\cr
                      {}\cr\noalign{\vskip0.6ex}
\left( \matrix{ \ a_{\sr_1}\ \cr \vdots \cr a_r\cr j}\right),\ \ j_{\sr C}=\infty & if $j>a_r$ or $C=\emptyset$\cr}$$
Put also  $\infty\ra C=C.$
The inductive extension of this operation 
 to a tableau
$T$  with $l$  columns for $j\in {\Na}^+\setminus \<T>$ given by
$$ j\Ra T=\left(j\ra T_1,
                  j_{\sr T_1}\Ra T_{2,l}\right) $$
is called the  insertion algorithm. 
\par
Note that the shape of $j\Ra T$ is the shape of $T$ obtained by adding one new corner. 
The entry of this corner is denoted by $j_{\sr T}.$ 
\par
This procedure (like many others used here)
is described in the wonderful book of B.E. Sagan ([Sa]). 
\parno
{\bf 2.6}\ \ \ 
Let $w=[a_{\sr 1},a_{\sr 2},\ldots,a_n]$ be a word.
According to Robinson-Schensted algorithm we associate an 
ordered pair
of tableaux $(T(w),Q(w))$ to $w.$ The procedure is 
fully explained 
in many places, for example, in [Sa, \S 3], [Kn, 5.1.4] or [F,4.1]. 
Here we explain only the inductive procedure
of constructing the first tableau $T(w)$ by insertions from the 
left. In what follows we call it RS procedure.
\item{(1)} \ \ \ Set $ {}_1T(w) = (a_n).$
\item{(2)} \ \ \ Set ${}_{j+1}T(w)=a_{n+1-j}\Ra{}_jT(w) \ .$
\item{(3)} \ \ \ Set $T(w)={}_nT(w) \ .$
\par
For example, let $w = [ 2 , 5,  1, 4, 3 ] $, then
$${}_1T(w) =
    \vcenter{
\halign{& \hfill#\hfill
\tabskip4pt\cr
\multispan{3}{\hrulefill}\cr
\ss1
\vb & 3  & \ts\vb\cr
\vs1
\multispan{3}{\hrulefill}\cr}} \quad
{}_2T(w) =
\vcenter{
\halign{& \hfill#\hfill
\tabskip4pt\cr
\multispan{3}{\hrulefill}\cr
\ss1
\vb &3 & \ts\vb\cr
\vs1
&&\cr
\ss1
\vb & 4 &  \ts\vb\cr
\vs1
\multispan{3}{\hrulefill}\cr}}\quad
{}_3T(w) =
\vcenter{
\halign{& \hfill#\hfill
\tabskip4pt\cr
\multispan{5}{\hrulefill}\cr
\ss1
\vb&1 & & 3  & \ts\vb\cr
\vs1
&&\multispan{3}{\hrulefill}\cr
\ss1
\vb & 4 &  \ts\vb\cr
\vs1
\multispan{3}{\hrulefill}\cr}}$$
$${}_4T(w) =
\vcenter{
\halign{& \hfill#\hfill
\tabskip4pt\cr
\multispan{5}{\hrulefill}\cr
\ss1
\vb & 1 && 3  & \ts\vb\cr
\vs1
&&\multispan{3}{\hrulefill}\cr
\ss1
\vb& 4 & \ts\vb\cr
\vs1
&&\cr
\ss1
\vb& 5 & \ts\vb\cr
\vs1
\multispan{3}{\hrulefill}\cr}}\qquad 
T(w)={}_5T(w)=
\vcenter{
\halign{& \hfill#\hfill
\tabskip4pt\cr
\multispan{5}{\hrulefill}\cr
\ss1
\vb&1 &&3  & \ts\vb\cr
\vs1
&&&&\multispan{0}{\hrulefill}\cr
\ss1
\vb& 2 & &4 & \ts\vb\cr
\vs1
&&\multispan{3}{\hrulefill}\cr
\ss1
\vb&5 & \ts\vb\cr
\vs1
\multispan{3}{\hrulefill}\cr}} $$
\par
The result due to Robinson and Schensted  implies the map
$\varphi: w \mapsto T(w)$ is a surjection from  ${{\Bbb S}}_n$ onto 
${\bR}_n$.
\parno
{\bf 2.7}\ \ \  For $T\in\bT_n$
one has (cf. for example, [M1, 2.4.14]) $\tau(T(w))=\tau(w).$
Thus, by 2.1 one has 
\proclaim Lemma. Let $S,T\in \bT_n.$ If $T\do S$ then
$\tau(T)\steq \tau(S).$
\par 
\parno
{\bf 2.8}\ \ \ 
Let us describe a few algorithms connected to RS
procedure which we use for proofs and constructions. 
\par
First let us describe some operations for columns and tableaux.
Consider a column $C=\pmat{a_{\sr 1}\cr\vdots\cr}.$
\item{(i)} For $m\in \<C>$ set $C\setminus \{m\}$ to be a column obtained from
 $C$ by deleting $m.$
\item{(ii)} For $j \in {\Na},\ j\not\in \<C>$ set $C+\{j\}$ to be a column obtained from $C$
by adding $j$ at the right place of $C$, that is if $a_i$ is the greatest
element of $\<C>$ smaller than $j$ then  $C+\{j\}$ is obtained from $C$ by adding $j$ between 
$a_i$ and $a_{i+1}.$ 
\item{(iii)} We define a pushing left operation. Again let 
$j \in {\Na},\ j\not\in \<C>$ and $j>a_{\sr 1}.$ 
Let $a_i$
be the greatest entry of $C$ smaller than $j$ and set :
$$C \la j:=\pmat{a_1\cr\vdots\cr a_{i-\sr 1}\cr j\cr a_{i+1}\cr\vdots\cr},\qquad j^{\sr C}:
                                                                     =a_i.$$
The last operation is extended to a tableau $T$ by induction
on the number of columns. Let $T_m$ be the last column of $T$ and assume $T_m^1<j.$ Then
$T\la j=( T_{1,m-1} \la j^{\sr T_m}, T_m \la j).$
We denote by $j^{\sr T}$ the element pushed out 
from the first column of the tableau in the last step.
\parno 
{\bf 2.9}\ \ \ 
The pushing left operation gives us a procedure of deleting a corner inverse to the
insertion algorithm. This is also described in many places, in particular, in all three 
books mentioned above.
\par
 As a result of insertion, we get a new tableau of a
shape obtained from the old one just by adding one corner. As a result of 
deletion, we  get a new tableau of a shape obtained from the old one
by removing one corner.
\par
Let $T=(T_1,\ldots, T_l).$ Recall the definition  of $\om_i(T)$ from 
2.3. 
Assume $\l_i>\l_{i+1}$ and let $c=c(\l_i,i)$ be a corner of $T$  on the $i$-th column.
To delete the  corner $c$ we
delete $\om_i(T)$ from the column $T_i$ and
push it left through the tableau $T_{1,i-1}.$ 
The element pushed out from
the tableau is denoted by $c^{\sr T}.$ 
This is written
$$T\La c:=\left(T_{1,i-1} \la \om_i(T),T_i\setminus\{\om_i(T)\},T_{i+1,l}\right)$$
For example,
$$\vcenter{
\halign{& \hfill#\hfill
\tabskip4pt\cr
\multispan{5}{\hrulefill}\cr
\ss1
\vb&1 &&3  & \ts\vb\cr
\vs1
&&&&\multispan{0}{\hrulefill}\cr
\ss1
\vb& 2 & &4 & \ts\vb\cr
\vs1
&&\multispan{3}{\hrulefill}\cr
\ss1
\vb&5 & \ts\vb\cr
\vs1
\multispan{3}{\hrulefill}\cr}}\La c(2,2)=
\vcenter{
\halign{& \hfill#\hfill
\tabskip4pt\cr
\multispan{5}{\hrulefill}\cr
\ss1
\vb & 1 && 3  & \ts\vb\cr
\vs1
&&\multispan{3}{\hrulefill}\cr
\ss1
\vb&4 & \ts\vb\cr
\vs1
&&\cr
\ss1
\vb&5 & \ts\vb\cr
\vs1
\multispan{3}{\hrulefill}\cr}},\qquad c^{\sr T}=2.$$ 
\par
Note that insertion and deletion are  indeed inverse since 
 for any $T\in {\Bbb T}_n$
$$c^{\sr T}\Ra (T\La c)=T\quad {\rm and} 
\quad(j\Ra T)\La j_{\sr T}=T\quad ({\rm for}\ j\not\in \<T>)$$
Note that sometimes we will write $T\La a$ where $a$ is a corner entry 
just as we have written above.
\par
Let $\{c_i\}_{i=1}^j$ be a set of corners of $T.$ By Robinson-Schensted procedure, one has  
$$\Cscr_T=\coprod_{i=1}^j\ \coprod_{y\pr\in\Cscr_{T\La c_i}}[c_i^{\sr T},y\pr]. \eqno{(*)}$$
\parno 
{\bf 2.10}\ \ \ Let us describe the jeu de taquin procedure (see [Sch]) which removes $T^i_j$
from $T.$ 
The resulting tableau is denoted by $T\setminus\{T^i_j\}.$ The idea of jeu de taquin is 
to remove $T^i_j$ from the tableau and to fill the gape created so that the resulting object
is again a tableau. The procedure goes as following. Remove a box from the tableau. Examine
the content of the box to the right of the removed box and that of the box below of the removed box.
Slide the box containing the smaller of these two numbers to the vacant position. Now repeat this
procedure to fill the hole created by the slide. Repeat the process until no holes remain, that is
until the hole has worked itself to the corner of the tableau.
\parno
The result due to M. P. Sch\"utzenberger [Sch] gives
\proclaim Theorem. If $T$ is a Young tableau then $T\setminus\{T^i_j\}$
is a Young tableau and the elimination of  different entries from $T$
by jeu de taquin is independent of the order chosen.
\par
Therefore, given $i_1,\ldots,i_s\in\<T>$, a tableau $T\setminus\{i_1,\ldots,i_s\}$ is a well defined 
tableau.

For example, let  us take
$$ T =
\vcenter{
\halign{& \hfill#\hfill
\tabskip4pt\cr
\multispan{7}{\hrulefill}\cr
\ss1
\vb & 1 &  &  2 & &  5 & \ts\vb\cr
\vs1
&&&&\multispan{3}{\hrulefill}\cr
\ss1
\vb & 3 &  &  4 & \ts\vb\cr
\vs1
&&\multispan{3}{\hrulefill}\cr
\ss1
\vb & 6 & \ts\vb\cr
\vs1
\multispan{3}{\hrulefill}\cr}}$$
Then a few tableaux obtained from $T$ by jeu de taquin are
$$T\setminus\{6\} =
\vcenter{
\halign{& \hfill#\hfill
\tabskip4pt\cr
\multispan{7}{\hrulefill}\cr
\ss1
\vb & 1 &  &  2 & &  5 & \ts\vb\cr
\vs1
&&&&\multispan{3}{\hrulefill}\cr
\ss1
\vb & 3 &  &  4 & \ts\vb\cr
\vs1
\multispan{5}{\hrulefill}\cr}},\quad
T\setminus\{3\} =
\vcenter{
\halign{& \hfill#\hfill
\tabskip4pt\cr
\multispan{7}{\hrulefill}\cr
\ss1
\vb & 1 &  &  2 & &  5 & \ts\vb\cr
\vs1
&&\multispan{5}{\hrulefill}\cr
\ss1
\vb & 4 & \ts\vb\cr
\vs1
&&\cr
\ss1
\vb & 6 & \ts\vb\cr
\vs1
\multispan{3}{\hrulefill}\cr}},\quad
T\setminus\{1,2\} =
\vcenter{
\halign{& \hfill#\hfill
\tabskip4pt\cr
\multispan{7}{\hrulefill}\cr
\ss1
\vb & 3 &  &  4 & &  5 & \ts\vb\cr
\vs1
&&\multispan{5}{\hrulefill}\cr
\ss1
\vb & 6 & \ts\vb\cr
\vs1
\multispan{3}{\hrulefill}\cr}}.$$
\parno 
{\bf 2.11}\ \ \ Given $s,t\ :\ 1\leq s<t\leq n$, set ${\Bbb T}_{\<s,t>}$
to be a set of Young tableaux with the entries $\{m_k\}_{k=s}^t$.
Let us define 
projection  $\pi_{s,t}:{\Bbb T}_n\ra {\Bbb T}_{\<s,t>}$
by  
$\pi_{s,t}(T)=T\setminus\{m_{\sr 1},\ldots,m_{s-1},m_{t+1},\ldots,m_n\}$. 
As a straightforward corollary of 2.10 (cf. for example, [M1, 4.1.1]), we get
\proclaim Theorem. for any $s,t\ :\ 1\leq s<t\leq n$ one has 
$\pi_{s,t}(T(w))=T(\pi_{s,t}(w)).$
\par
\parno
{\bf 2.12}\ \ \ As a  straightforward corollary of lemma 2.2 (iv) and theorem 2.11, we get that 
$\do$ is preserved 
under projections and, as a straightforward corollary of lemma 2.2 (i) and RS procedure,
we get that $\do$ is preserved 
under insertions, namely
\proclaim Proposition. Let $T,S$ be in ${\Bbb T}_n.$ If $T\do S$ then
\item {(i)} for any $s,t\ :\ 1\leq s<t\leq n$ one has $\pi_{s,t}(T)\do \pi_{s,t}(S).$
\item {(ii)} for any $a\not\in \{m_s\}_{s=1}^n$ one has $a\Ra T\do a\Ra S.$  
\par
\parno
{\bf 2.13}\ \ \ Consider $T\in \bT_n.$
Note that
$$\pi_{i,i+1}(T)=\cases{
{\vcenter{
\halign{& \hfill#\hfill
\tabskip4pt\cr
\multispan{3}{\hrulefill}\cr
\ss1
\vb & {\it i} & \ts\vb\cr
\vs1
&&\cr
\ss1
\vb & {\it i}+1 & \ts\vb\cr
\vs1
\multispan{3}{\hrulefill}\cr}}} & if $i\in\tau(T)$\cr
&\cr
\noalign{\medskip}\cr
{\vcenter{
\halign{& \hfill#\hfill
\tabskip4pt\cr
\multispan{5}{\hrulefill}\cr
\ss1
\vb & {\it i} && {\it i}+1 & \ts\vb\cr
\vs1
\multispan{5}{\hrulefill}\cr}}} & if $i\not\in\tau(T)$\cr}$$
We need the following properties of the chain order.
\proclaim Proposition. Let $S,T\in \bT_n.$ 
\item {(i)} If $T\cho S$ then $\tau(T)\st\tau(S).$
\item {(ii)}If $T\cho S$ then for any $i,j\ :\ 1\leq i<j\leq n$ one has $\pi_{i,j}(T)\cho \pi_{i,j}(S).$ 
\item {(iii)} If $T\do S$ then $T\cho S.$
\par
\Pf
The first two assertions are trivial. The third assertion is a corollary of 
Steinberg's construction explained in 1.9 and 1.10 and of proposition 2.12 (i) or of
the results explained in 1.11.

Indeed, $y\do w$ implies that $\ov\Oscr_y\supseteq\ov\Oscr_w.$ 
Thus, $T\do S$ implies $\sh(T)\leq \sh(S).$ 
By proposition 2.12 (i),  $T\do S$ implies 
$\pi_{s,t}(T)\do \pi_{s,t}(S)$ for any $s,t\ :\ 1\leq s<t\leq n.$ Altogether, 
this provides $T\cho S.$ 
\QED
\parno
{\bf 2.14}\ \ \ All the results for the tableaux with two columns
can be translated to tableaux with two rows by Schensted-Sch\"utzenberger
theorem (cf. [Kn, 5.4.1]).
\proclaim Theorem. For any $w\in {\Bbb S}_n$
one has $T^\dagger(w)=T(\ov w).$
\par
\bigskip
\bigskip
{\bf 3. Combinatorics of $\bT_n^2.$} 
\bigskip
\parno
{\bf 3.1}\ \ \ 
Recall from 1.5 that  ${\Bbb T}_n^2\st {\Bbb T}_n$  is the set of 
Young  tableaux
with 2 columns. For $T\in {\Bbb T}_n^2$ let $\l_1(T)$ be the length of the first column
and $\l_2(T)$ be the length of the second column, that is $\sh(T)=(\l_1(T),\l_2(T)).$ 
\proclaim Lemma. Let $T\in \bT_n^2$ be such that  $c_{\sr T}(n)=2.$ Set $T\pr=\pi_{\sr 2,n}(T).$
Then $c_{\sr T\pr}(n)=2$ if and only if either $\l_1 (T)>\l_2(T)$ or there exists $i$
such that $T_2^i<T_1^{i+1}.$
\par
The proof is a straightforward and easy computation, so we omit it.
\parno
{\bf 3.2}\ \ \ 
Consider tableaux $T, S\in \bT_n^2.$
\proclaim Lemma. If $c_{\sr T}(n)=1$ and $S\chgs T$  then $c_{\sr S}(n)=1.$
\par
\Pf
This is true for $n=3.$ Assume this is true for $k=n-1$ and show for $k=n.$
If $c_{\sr T}(n)=1$ then $\l_1(T)>\l_2(T).$ Since $S\chgs T$ one has $\sh (S)>\sh (T).$
Thus, $\l_1(S)>\l_2(S).$
Assume $c_{\sr S}(n)=2.$ Then by lemma 3.1 $c_{\sr \pi_{2,n}(S)}(n)=2.$  On the other hand,
$c_{\sr \pi_{2,n}(T)}(n)=1$ by the induction assumption, and this is a contradiction.
\QED
\parno
{\bf 3.3}\ \ \ As a corollary of lemma 3.2, we get
\proclaim Corollary.  For $S,T\in\bT_n^2$ one has
\item{(i)} If $T\ne S$ and $\sh T=\sh S$  
then $T$ and $S$ are incompatible in the chain order.
\item{(ii)} If $S\chgs T$ then $\<S_1>\supset \<T_1>$ and $\<S_2>\subset \<T_2>.$
\par
\Pf
(i) This is true for $n=3.$ 
Assume this is true for $n-1$ and show for $n.$
\item {(a)} If $c_{\sr T}(n)=c_{\sr S}(n)$ then
$\pi_{1,n-1}(T)\ne\pi_{1,n-1}(S)$ and  $\sh \pi_{1,n-1}(T)=
\sh \pi_{1,n-1}(S)$, hence, they are incompatible by assumption hypothesis.
\item {(b)} If $c_{\sr T}(n)=1$ and $c_{\sr S}(n)=2$ then
  $\sh \pi_{1,n-1}(T)=(\l_1(T)-1,\l_2(T))$ and\hfil\break
$\sh \pi_{1,n-1}(S)=(\l_1(T),\l_2(T)-1)$ so that $\sh \pi_{1,n-1}(T)<\sh \pi_{1,n-1}(S).$
Hence, $S\not \cho T.$
On the other hand, by lemma 3.2 $T\not\cho S.$

(ii)\ \ For any $j\ :\ j<n$ one has $\pi_{1,j}(S)\chg \pi_{1,j}(T)$. 
If $c_{\sr T}(j)=1$ then by lemma 3.2 applied to $\pi_{1,j}(T),\ \pi_{1,j}(S)$ 
we get $c_{\sr S}(j)=1.$ Further note that $\<T_2>=\{i\}_{i=1}^n\setminus \<T_1>.$ 
\QED
Note that in general neither of these assertion is true, as it is shown in the 
following example: $T\chos S$ where
$$T=\vcenter{
\halign{& \hfill#\hfill
\tabskip4pt\cr
\multispan{7}{\hrulefill}\cr
\ss1
\vb&1 &&2  && 5 & \ts\vb\cr
\vs1
&&&&\multispan{3}{\hrulefill}\cr
\ss1
\vb& 3 & &4 & \ts\vb\cr
\vs1
&&\multispan{3}{\hrulefill}\cr
\ss1
\vb&6 & \ts\vb\cr
\vs1
\multispan{3}{\hrulefill}\cr}}\qquad {\rm and}\qquad 
S=\vcenter{
\halign{& \hfill#\hfill
\tabskip4pt\cr
\multispan{7}{\hrulefill}\cr
\ss1
\vb&1 &&2  && 5 & \ts\vb\cr
\vs1
&&&&\multispan{3}{\hrulefill}\cr
\ss1
\vb& 3 & &6 & \ts\vb\cr
\vs1
&&\multispan{3}{\hrulefill}\cr
\ss1
\vb&4 & \ts\vb\cr
\vs1
\multispan{3}{\hrulefill}\cr}} $$
Therefore, to avoid two tableaux of the same shape to be in
the chain order we have to restrict the chain order
by the demand that if for some $T\chos S$ and for some $i,j\ :\ 1\leq i<j\leq n$
one has $D_{\<i,j>}(T)=D_{\<i,j>}(S)$ then $\pi_{i,j}(T)=\pi_{i,j}(S).$ As we see,
we do not need this restriction on $\bT_n^2.$
\parno
{\bf 3.4}\ \ \ One has
\proclaim Lemma. If $T\dos S$ and $c_{\sr T}(n)=2,\ c_{\sr S}(n)=1$ then $T\pr\do S$
where $T\pr=(T_1+\{n\},T_2\setminus\{n\}).$
\par
\Pf
Indeed, if $T\dos S$ then by proposition 2.12 (i) $\pi_{1,n-1}(T)\do \pi_{1,n-1}(S)$ and further by 
proposition 2.12 (ii) $(T_1+\{n\},T_2\setminus\{n\})=n\Ra \pi_{1,n-1}(T)\do n\Ra \pi_{1,n-1}(S)=S.$
\QED
\parno
{\bf 3.5}\ \ \ Now we construct the special representative of $\Cscr_T$ which plays the key 
role in our constructions.
 
Given $T\in\bT_n^2$, put $T_{(n)}=T.$ Let $z_i:=\max\<T_{(i)}>.$ Obviously $z_i$ is a corner element
of $T_{(i)}.$ Set $T_{(i-1)}=T_{(i)}\La z_i.$
Recall notion $c^{\sr T}$ from 2.9. 
Set $a_i:=z_i^{\sr T}.$ 
\par
Note that for any $s\in T_2$ there exists a unique $i$
such that $s=z_i=z_{i-1}.$ For $s\in T_2$
set  $T\{s\}:=T_{(i)}.$ 
\par
For example, let
$$T=\vcenter{
\halign{& \hfill#\hfill
\tabskip4pt\cr
\multispan{5}{\hrulefill}\cr
\ss1
\vb & 1 &  &  3 &\ts\vb\cr\vs1
&&&&\cr
\ss1
\vb & 2 &  &  5 & \ts\vb\cr\vs1
&&&&\cr
\ss1
\vb & 4 &  &  6 & \ts\vb\cr\vs1
&&\multispan{3}{\hrulefill}\cr
\ss1
\vb & 7 &\ts\vb\cr\vs1
\multispan{3}{\hrulefill}\cr}}$$
then
$$T_{(7)}=T,\ z_7=7;\quad T_{(6)}=T_{(7)}\La 7 =\vcenter{
\halign{& \hfill#\hfill
\tabskip4pt\cr
\multispan{5}{\hrulefill}\cr
\ss1
\vb & 1 &  &  3 &\ts\vb\cr\vs1
&&&&\cr
\ss1
\vb & 2 &  &  5 & \ts\vb\cr\vs1
&&&&\cr
\ss1
\vb & 4 && 6& \ts\vb\cr\vs1
\multispan{5}{\hrulefill}\cr}},\ a_7=7,\ z_6=6,\ T\{6\}=T_{(6)};$$

$$
T_{(5)}=T_{(6)}\La 6=\vcenter{
\halign{& \hfill#\hfill
\tabskip4pt\cr
\multispan{5}{\hrulefill}\cr
\ss1
\vb & 1 &  &  3 &\ts\vb\cr\vs1
&&&&\cr
\ss1
\vb & 2 &  &  5 & \ts\vb\cr\vs1
&&\multispan{3}{\hrulefill}\cr
\ss1
\vb & 6 & \ts\vb\cr\vs1
\multispan{3}{\hrulefill}\cr}},\ a_6=4,\ z_5=6; $$

$$T_{(4)}=T_{(5)}\La 6=\vcenter{
\halign{& \hfill#\hfill
\tabskip4pt\cr
\multispan{5}{\hrulefill}\cr
\ss1
\vb & 1 &  &  3 &\ts\vb\cr\vs1
&&&&\cr
\ss1
\vb & 2 && 5 & \ts\vb\cr\vs1
\multispan{5}{\hrulefill}\cr}},\ a_5=6,\ z_4=5,\ T\{5\}=T_{(4)};$$

$$T_{(3)}=T_{(4)}\La 5=\vcenter{
\halign{& \hfill#\hfill
\tabskip4pt\cr
\multispan{5}{\hrulefill}\cr
\ss1
\vb & 1 &  &  3 &\ts\vb\cr\vs1
&&\multispan{3}{\hrulefill}\cr
\ss1
\vb & 5 & \ts\vb\cr\vs1
\multispan{3}{\hrulefill}\cr}},\ a_4=2,\ z_3=5;$$

$$T_{(2)}=T_{(3)}\La 5=\vcenter{
\halign{& \hfill#\hfill
\tabskip4pt\cr
\multispan{5}{\hrulefill}\cr
\ss1
\vb & 1 & &  3 & \ts\vb\cr\vs1
\multispan{5}{\hrulefill}\cr}},\ a_3=5, \ z_2=3,\ T\{3\}=T_{(2)};$$

$$
 T_{(1)}=T_{(2)}\La 3=\vcenter{
\halign{& \hfill#\hfill
\tabskip4pt\cr
\multispan{3}{\hrulefill}\cr
\ss1
\vb &  3 &\ts\vb\cr\vs1
\multispan{3}{\hrulefill}\cr}},\ a_2=1,\ z_1=a_1=3.$$
Put $w_{\sr T}:=[a_n,a_{n-1},\ldots,a_1].$ In our example 
$w_{\sr T}=[7,4,6,2,5,1,3].$
Note that by Robinson-Schensted procedure $T(w_{\sr T})=T.$ 
\par
Now we can formulate the main theorem of the paper
\proclaim Theorem. For $T,S\in\bT_n^2$ one has $T\chos S$ iff $w_{\sr T}\dos w_{\sr S}.$
\par
To prove the theorem we need a few technical lemmas. 

\par
\parno
{\bf 3.6}\ \ \ First of all we show that $w_{\sr T}$ is a maximal element of
$\Cscr_T$ in the weak order.
\proclaim Lemma. For any $y\in\Cscr_T$ one has $y\do w_{\sr T}.$
\par
\Pf
This is true for $n=3.$ Assume that this is true for $k\leq n-1$ and show for $n.$
Take $T\in\bT_n^2.$ Set $\om_{\sr 1}:=\om_{\sr 1}(T)$ and $\om_{\sr 2}:=\om_{\sr 2}(T).$
\item {(i)} If $c_{\sr T}(n)=1$ (which means  $\om_{\sr 1}=n$) then $w_{\sr T}=[n,w_{\sr \pi_{1,n-1}(T)}]$ and for any
$y$ such that $T(\pi_{1,n-1}(y))=\pi_{1,n-1}(T)$ one 
has by lemma 2.2 (ii) $y\do [n,\pi_{1,n-1}(y)].$
In particular, for any $y\in \Cscr_T$ one has $y\do [n,\pi_{1,n-1}(y)]\do 
[n,w_{\sr \pi_{1,n-1}(T)}]$ just by induction assumption and lemma 2.2 (i). 
\item {(ii)} If $c_{\sr T}(n)=2$ (which means $\om_{\sr 2}=n$) then $\om_{\sr 1}^{\sr T}=\om_{\sr 2}^{\sr T}=\om_{\sr 1}$
thus, by 2.9 $(*)$ any $y\in \Cscr_T$ has a form $y=[\om_{\sr 1},y\pr]$
where either $T(y\pr)=T\La n=:T\pr$ or $T(y\pr)=T\La\om_{\sr 1}=:T\prpr.$
Note that 
$w_{\sr T}=[\om_{\sr 1},w_{\sr T\pr}]$ thus, by induction assumption
and lemma 2.2 (i)
for any $y\pr\in\Cscr_{T\pr}$ one has $[\om_{\sr 1},y\pr]\do w_{\sr T}.$
For any $y\pr\in\Cscr_{T\prpr}$ one has just by induction assumption
that $y\pr\do w_{\sr T\prpr}$ where $w_{\sr T\prpr}=[\om_{\sr 1}(T\prpr),n,z],$ where
$z$ is the rest of this word.
Note that by definition of the right weak order $w_{\sr T\prpr}\dos [n,\om_{\sr 1}(T\prpr),z]$ so that for any 
$y\pr\in\Cscr(T\prpr)$ one has $y\pr\dos [n,\om_{\sr 1}(T\prpr),z].$ On the other hand,
$T([n,\om_{\sr 1}(T\prpr),z])=T\pr.$ Indeed, $T(\pi_{1,n-2}([n,\om_{\sr 1}(T\prpr),z]))=\pi_{1,n-2}(T\prpr)=
\pi_{1,n-2}(T\pr)$ and by RS procedure  $c_{\sr T([n,\om_{\sr 1}(T\prpr),z])}(n)=1.$
Thus, for any $y\pr\in\Cscr(T\prpr)$ one has $y\pr\dos [n,\om_{\sr 1}(T\prpr),z]\do w_{\sr T\pr}.$
Applying lemma 2.2 (i) we get that for any 
$y\pr\in \Cscr(T\prpr)$ one has $[\om_{\sr 1},y\pr]\dos w_{\sr T}.$
\QED
\parno
{\bf 3.7}\ \ \ As a corollary of lemma 3.6 and its proof, we get 
\proclaim Corollary. If $c_{\sr T}(n)=2$ and $T\pr=n\Ra \pi_{1,n-1}(T)$  
then for any $y\in \Cscr_T$ one has that $y\dos w_{\sr T\pr}.$
\par
\Pf
Indeed, $y\do w_{\sr T}=[\om_{\sr 1}(T),n,z]\dos [n,\om_{\sr 1}(T), z]$ and as we have shown in
(ii) of the proof of lemma 3.6, $[n,\om_{\sr 1}(T), z]\in\Cscr(T\pr)$, hence, by lemma 3.6,
$y\dos w_{\sr T\pr}.$
\QED
\parno
{\bf 3.8}\ \ \ 
Let us return to the description of the orders on the level of tableaux. 
\proclaim Lemma. Let $T,S\in \bT_n^2.$ If $S\chgs T$  and $c_{\sr T}(n)=c_{\sr S}(n)$
then $\om_{\sr 1}(T)=\om_{\sr 1}(S).$
\par
\Pf
If $c_{\sr T}(n)=1$ then $\om_{\sr 1}(T)=\om_{\sr 1}(S)=n.$ 
\par
Assume that  $c_{\sr T}(n)=2.$
For $n=4$ this is true. Assume this is true for $k=n-1$ and show for $k=n.$
Consider $T\pr=\pi_{1,n-1}(T)$ and $S\pr=\pi_{1,n-1}(S).$ By proposition 2.13 (ii),
 $S\pr\chgs T\pr.$
\item {(i)} If $c_{\sr S\pr}(n-1)=2$ then by lemma 3.2 $c_{\sr T\pr}(n-1)=2$
and by induction assumption $\om_{\sr 1}(S\pr)=\om_{\sr 1}(T\pr).$ On the other hand,
$\om_{\sr 1}(S\pr)=\om_{\sr 1}(S)$ and $\om_{\sr 1}(T\pr)=\om_{\sr 1}(T).$
\item {(ii)} If $c_{\sr S\pr}(n-1)=1$ then $n-1\not\in\tau (S).$ Thus, by proposition 2.13 (i)
$n-1\not\in\tau (T)$ so that $c_{\sr T\pr}(n-1)=1$ and $\om_{\sr 1}(T)=\om_{\sr 1}(S)=n-1.$
\QED
\parno
{\bf 3.9}\ \ \ Let $S$ be a tableau with two columns. 
For $x\in \<S_2>$ recall notion $S\{x\}$ from 3.5.
Since $x$ is $\om_{\sr 2}(S\{x\})$ we consider $S\{x\}\La x$ and
get $x^{\sr{S\{x\}}}$ (as defined in 2.9).
Obviously, $x^{\sr{S\{x\}}}$  is some element of $S_1.$  
\proclaim Lemma. Let $T,S\in \bT_n^2.$ $w_{\sr T}\dos w_{\sr S}$ iff 
$\<S_2>\st \<T_2>$ and for any $x\in \<S_2>$ one has
$x^{\sr S\{x\}}\in \<T_1>.$
\par
\Pf
First of all note that $w_{\sr T}\dos w_{\sr S}$ implies that $\<S_2>\st \<T_2>$ by corollary
3.3 (ii). As well, this implies that for any $x\in \<S_2>$ one has
$x^{\sr S\{x\}}\in \<T_1>.$ Indeed,  assume that there exist $x\in \<S_2>$ such that
$a:=x^{S\{x\}}\not\in \<T_1>.$
Then $c_{\sr T}(a)=2$ and by definition of $w_{\sr T}$
one has $p_{w_T}(a)>p_{w_T}(x).$ On the other hand, $p_{w_S}(a)=p_{w_S}(x)-1.$ 
Thus, we have found $a<x$ such that $p_{w_T}(x)<p_{w_T}(a)$ and $p_{w_S}(x)>p_{w_S}(a).$
This implies that $w_{\sr T}\not\do w_{\sr S}.$
\par
We show the other direction  by induction. The claim is true for $n=4.$ Assume this is true for 
$k\leq n-1$ and show for $k=n.$
\item {(i)} If $c_{\sr T}(n)=1$ then $c_{\sr S}(n)=1$ since $\<S_2>\st\<T_2>.$ 
Set $T\pr:=\pi_{1,n-1}(T)$ and $S\pr:=\pi_{1,n-1}(S).$  One has $T_1\pr=T_1-n$
and $S_2\pr=S_2.$ As well,  $x^{\sr S\{x\}}=x^{\sr S\pr\{x\}}$ for any $x\in\<S_2>.$
Thus, if $x^{\sr S\{x\}}\in \<T_1>$ then $x^{\sr S\pr\{x\}}\in \<T_1\pr>.$ By induction hypothesis,
this provides $w_{\sr T\pr}\dos w_{\sr S\pr}.$ Note that
$w_{\sr T}=[n,w_{\sr T\pr}],\ w_{\sr S}=[n,w_{\sr S\pr}].$
Thus, by lemma 2.2 (i) $w_{\sr T}\dos w_{\sr S}.$

\item {(ii)} Assume that $c_{\sr T}(n)=c_{\sr S}(n)=2.$  Since $\<S_2>\st \<T_2>,$ we get that 
$\<S_1>\supset \<T_1>$ and, in particular, $\om_{\sr 1}(T)\leq \om_{\sr 1}(S).$ 
Since $n^S=\om_{\sr 1}(S),$  by the condition $n^S\in \<T_1>$ we get that
$\om_{\sr 1}(T)=\om_{\sr 1}(S).$ Let us denote it by $\om_{\sr 1}.$ Thus, by the construction
 $w_{\sr T}=[\om_{\sr 1},n, w_{\sr T\pr}]$ and $w_{\sr S}=[\om_{\sr 1},n, w_{\sr S\pr}]$
where $T\pr=(T_1-\om_{\sr 1}, T_2-n)$ and $S\pr=(S_1-\om_{\sr 1}, S_2-n).$ 
Let us show that $T\pr, S\pr$ satisfy the conditions. It is obvious that $\<S\pr_2>\st \<T\pr_2>.$ 
Further, 
 $n^{\sr S}=\om_{\sr 1}\in\<T_1>$
and for any $x\ : x\ne n,\ x\in\<S_2>$ one has $x\in\<S\pr_2>$ and $x^{\sr S\{x\}}=x^{\sr S\pr\{x\}}$ just by
construction. Moreover,  for such $x$ one has $x^{\sr S\{x\}}\ne\om_{\sr 1}.$ Thus, the condition
$ x^{\sr S\{x\}}\in\<T_1>$ for any $x\in\<S_2>$
provides $x^{\sr S\pr\{x\}}\in\<T\pr_1>$ for any $x\in\<S\pr_2>.$ By induction hypothesis, this implies $w_{\sr T\pr}\dos w_{\sr S\pr}.$
Again by lemma 2.2 (i) if $w_{\sr T\pr}\dos w_{\sr S\pr}$ then $w_{\sr T}\dos w_{\sr S}.$ 

\item{(iii)} Finally, assume that $c_{\sr T}(n)=2$ and $c_{\sr S}(n)=1.$ 
Consider $T\pr=n\Ra\pi_{1,n-1}(T).$ Note that $\<S_2>\st \<T_2>$ and $n\not\in \<S_2>$
imply that $\<S_2>\st \<T\pr_2>.$ Let us show that $T\pr, S$ satisfy the second condition as well.
Indeed, $T\pr_1=(T_1+n).$ Thus,  
for any $x\in \<S_2>$ one has
$x^{S\{x\}}\in \<T_1>$ iff $x^{S\{x\}}\in \<T\pr_1>.$ By (i), this implies $w_{\sr S}\dgs w_{\sr T\pr}$
and by corollary 3.7 $w_{\sr T\pr}\dgs w_{\sr T}.$ This completes the proof.
\QED
\parno
{\bf 3.10}\ \ \ We need the following result
about the chain order 
\proclaim Lemma. Let $T,S\in \bT_n.$ Let $T\chos S$ and assume that $c_{\sr T}(n)=c_{\sr S}(n)=2.$
If $T\pr=T\La n,\ S\pr=S\La n$ then $T\pr\chos S\pr.$
\par
\Pf
By lemma 3.8, the assumption  $c_{\sr T}(n)=c_{\sr S}(n)=2$ implies
$\om_{\sr 1}(T)=\om_{\sr 1}(S)$ and we will denote it by $\om_{\sr 1}.$
We give a proof by induction.
This is true for $n=4.$ Assume this is true for $k=n-1$ and show for $k=n.$ 
\item {(i)} Suppose that $\om_{\sr 1}=n-1$ then $T\pr$ is equivalent
to $\pi_{1,n-1}(T)$ and $S\pr$ is equivalent to $\pi_{1,n-1}(S)$ and the
statement is obvious.
\item {(ii)} Let us consider the case $\om_{\sr 1}<n-1.$
We have that $T\pr,S\pr\in {\Bbb T}_{n-1}^2.$ To show that 
$T\pr\chos S\pr$ we note first that $\sh(T\pr)=(\lambda_1(T),\lambda_2(T)-1)$
and $\sh(S\pr)=(\lambda_1(S),\lambda_2(S)-1).$ Thus, $\sh(T\pr)<\sh(S\pr).$
As well, $\sh(\pi_{1,n-2}(T\pr))=(\lambda_1(T)-1,\lambda_2(T)-1)$ and 
$\sh(\pi_{1,n-2}(S\pr))=(\lambda_1(S)-1,\lambda_2(S)-1).$ Thus, again, $\sh(\pi_{1,n-2}(T\pr))<
\sh(\pi_{1,n-2}(S\pr)).$ Let us show that $\pi_{1,n-3}(T\pr)\chos \pi_{1,n-3}(S\pr)$ by induction
hypothesis. Indeed, $P=T,S$ and for $P\pr=T\pr, S\pr$
one has  $\pi_{1,n-3}(P\pr)=\pi_{1,n-3}(\pi_{1,n-1}(P)\La n-1)$
so that by induction hypothesis $\pi_{1,n-3}(T\pr)\chos\pi_{1,n-3}(S\pr).$
To complete the proof we have to show that $\pi_{2,n-1}(T\pr)\cho \pi_{2,n-1}(S\pr).$ 
Indeed, set $P\prpr=\pi_{2,n}(P)$ 
where $P$ is $T$ or $S$. Since $S\chgs T$ one has $\lambda_1(S)>\lambda_2(S)$. Thus, $S$
satisfies conditions (i) and (ii) of lemma 3.1, so that $c_{\sr S\prpr}(n)=2.$ 
This implies in turn by lemma 3.2 that $c_{\sr T\prpr}(n)=2.$ 
In particular, this provides  $\pi_{2,n-1}(P\pr)=P\prpr\La n.$ Hence, 
$\pi_{2,n-1}(T\pr)\cho \pi_{2,n-1}(S\pr)$
by induction assumption.
\QED
Note that this property is unique for $\bT_n^2$. Indeed, in general, the facts $T\chos S$ and $c_{\sr T}(n)=c_{\sr S}(n)$
even do not provide that $\<T\La n>=\<S\La n>.$ 
\parno 
{\bf 3.11}\ \ \ 
Now we are ready to prove theorem 3.5. Let us recall its formulation.
\proclaim Theorem. For $T,S\in\bT_n^2$ one has $T\chos S$ iff $w_{\sr T}\dos w_{\sr S}.$
\par
\Pf
As we explained in 1.11, $w_{\sr T}\dos w_{\sr S}$ implies $T\chos S.$
\par
We will show the other direction by induction.
For $n=3$ the other direction is true. Assume that for $k\leq n-1$ if $T\chos S$
then $w_{\sr T}\dos w_{\sr S}$ and show this for $k=n.$
\par
Assume $T\chos S.$ 
\item {(i)} If $c_{\sr S}(n)=1$ then consider $S\pr=\pi_{1,n-1}(S)$ and
$T\pr=\pi_{1,n-1}(T).$ By proposition 2.13 (ii), $T\pr\cho S\pr,$ thus, by 
induction assumption $w_{\sr T\pr}\do w_{\sr S\pr}.$
One has  $w_{\sr S}=[n,w_{\sr S\pr}]$
and by lemma 2.2 (i) this implies $[n,w_{\sr S\pr}]\dg [n,w_{\sr T\pr}]=w_{\sr T\prpr}$
where $T\prpr=n\Ra T\pr.$ By corollary 3.7, $w_{\sr T\prpr}\dg w_{\sr T}.$
Thus, we get in that case $w_{\sr T}\do w_{\sr S}.$
\item {(ii)} If $c_{\sr S}(n)=2$ then by lemma 3.2 $c_{\sr T}(n)=2.$ By lemma 3.8,
$\om_{\sr 1}(T)=\om_{\sr 1}(S)=:\om_{\sr 1}.$ As well, $\om_{\sr 1}=n^{\sr S}=n^{\sr T}.$
Consider $S\pr=S\La n$
and $T\pr=T\La n.$  By lemma 3.10, $T\pr\chos S\pr$ and by induction hypothesis this provides 
$w_{\sr T\pr}\dos w_{\sr S\pr}.$
On the other hand, $w_{\sr T}=[\om_{\sr 1},w_{\sr T\pr}]$ and $w_{\sr S}=[\om_{\sr 1}, w_{\sr S\pr}].$
Thus, by lemma 2.2 (i) in that case, as well, $w_{\sr T}\dos w_{\sr S}.$
\QED
\parno
{\bf 3.12}\ \ \ The first and very easy corollary of the theorem is
\proclaim Corollary. For $T,S\in \bT_n^2$ one has $T\dos S$ iff $w_{\sr T}\dos w_{\sr S}.$
\par
\Pf
The implication $w_{\sr T}\dos w_{\sr S}\ \Rightarrow\ T\dos S$ follow just from the definition; 
the other implication is obvious
from theorem 3.11 since $T\dos S$ implies $T\chos S.$
\QED
\parno
{\bf 3.13}\ \ \ 
As well, we get the following geometric fact from this purely combinatorial theorem 
\proclaim Corollary. For orbital varieties $\Vscr_T, \Vscr_S$ of nilpotent order 2 one has $\Vscr_S\st\ov\Vscr_T$
if and only if $\gn\cap^{w_{\sr S}}\gn\st \gn\cap^{w_{\sr T}}\gn$ that 
the inclusion of orbital variety closures is determined by inclusion
of generating subspaces.
\par
\Pf
Again one implication is obvious from the definition and the other one from theorem 3.11 
since $\Vscr_S\st\ov\Vscr_T$ implies by 1.11 $S\chgs T.$ 
\QED
\parno
{\bf 3.14}\ \ \ Theorem 3.11 and corollary 3.12 provide us also
\par
\proclaim Corollary. $\cho$ and $\do$ coincide on orbital varieties of nilpotent order 2.
\par
\parno
{\bf 3.15}\ \ \ 
Note that lemma 3.9  together with theorem 3.11 give the exact description of inclusion of
orbital variety closures of nilpotent order $2$ in terms of Young tableaux. 
Since $\go,\, \do,$ and $\cho$ coincide on $\bT_n^2$ we will denote them simply by $\leq$
and the cover in $\leq$ simply by $\Dscr(T).$ 

Let us first give the recursive description of $S\ :\ S>T$ and of $\Dscr(T)$ for $T\in \bT_n^2.$
\proclaim Proposition. {Let $T\in \bT_n^2.$ One has 
\item {(i)} If $c_{\sr T}(n)=1$ then $S>T$ iff $S=n\Ra S\pr$ where $S\pr>\pi_{1,n-1}(T).$
In particular, $\Dscr(T)=\{n\Ra S\}_{S\in\Dscr(\pi_{1,n-1}(T))}.$
\item {(ii)} If $c_{\sr T}(n)=2$ then $S>T$ in the next two cases.
Either $S=n\Ra S\pr$ where $S\pr\geq\pi_{1,n-1}(T)$
or $S=\om_{\sr 1}(T)\Ra S\pr$ where $S\pr>T\La n.$ In particular, 
$\Dscr(T)=\{\om_{\sr 1}(T)\Ra(n\Ra S\prpr)\}_{S\prpr\in \Dscr(T_1\setminus\{om_1(T)\},T_2\setminus\{n\})}
\cup \{(T_1+\{n\},T_2\setminus\{n\})\}.$\hbox{}
\item{} In particular, for any $T\in \bT_n^2$ and for any $S\in\Dscr(T)$ one has 
$\sh(S)=(\l_1(T)+1,\l_2(T)-1).$}
\par
\Pf
Indeed, if $c_{\sr T}(n)=1$ and $S>T$ then by lemma 3.2 one has $c_{\sr S}(n)=1.$ Thus,
$w_{\sr T}=[n,w_{\sr \pi_{1,n-1}(T)}]$ and $w_{\sr S}=[n,w_{\sr \pi_{1,n-1}(S)}].$
One has by lemma 2.2 (i) that $w_{\sr T}\dos w_{\sr S}$ iff 
$w_{\sr \pi_{1,n-1}(T)}\dos w_{\sr \pi_{1,n-1}(S)}$ which is equivalent by theorem 3.11 and 
its corollaries to $\pi_{1,n-1}(T)<\pi_{1,n-1}(S).$ 
Now if $c_{\sr T}(n)=1$ then $\sh (\pi_{1,n-1}(T))=(\l_1(T)-1,\l_2(T))$ and for any
 $S\in\Dscr(T)$ one has $\sh (\pi_{1,n-1}(S))=(\l_1(S)-1,\l_2(S)).$
Note that $\pi_{1,n-1}(S)>\pi_{1,n-1}(T)$ by shape consideration. The same shape considerations
show that if $S\in\Dscr(T)$ then $\pi_{1,n-1}(S)\in\Dscr(\pi_{1,n-1}(T))$ and that
for any $S\pr\in \Dscr(\pi_{1,n-1}(T))$ one has $n\Ra S\pr\in \Dscr(T).$
\par
Now assume that $c_{\sr T}(n)=2.$ Consider $S\ :\ S>T.$ If $c_{\sr S}(n)=2$ then by
lemma 3.10 and corollary 3.14 $(T\La n)<(S\La n).$ Thus, $S=\om_{\sr 1}(T)\Ra S\pr$ where $S\pr>T\La n.$
If $c_{\sr S}(n)=1$ then by lemma 3.4 and corollary 3.14 $S\geq T\pr=(T_1+\{n\},T_2\setminus\{n\}).$ 
Thus, by (i)
$S=n\Ra S\pr$ where $S\pr\geq\pi_{1,n-1}(T).$
If $S\in \Dscr(T)$ then
\item {(a)} If $c_{\sr S}(n)=1$ then by lemma 3.4  $S=(T_1+\{n\},T_2\setminus\{n\})$
\item {(b)} If $c_{\sr S}(n)=2$ then by lemma 3.10 and (i) $S=\om_{\sr 1}(T)\Ra(n\Ra S\prpr)$
where $S\prpr\in \Dscr(T_1\setminus \{\om_1(T)\},T_2\setminus \{n\}).$
\par
The note on the shape of $S\in\Dscr(T)$ is obvious.
\QED
{\bf 3.16}\ \ \ Let us give explicit description of $\Dscr(T).$
Consider $T\in\bT_n^2.$ One can write $T_2$ as the union of connected subsequences 
$\<T_2>=\{a_1,a_1+1,\ldots a_1+k_1\}\cup\ldots\cup \{a_s,\ldots,a_s+k_s\}$
where $a_i>a_{i-1}+k_{i-1}+1$ for any $i\ :\ 1<i\leq s.$
For any $x\in \<T_2>$ set $T\<x>:=(T_1+\{x\},T_2\setminus\{x\}).$ Note that $T\<x>$ is always a 
tableau. Recall notion of $T\{x\}$ from 3.5. Note that for $x\in \<T_2>$ sometimes $T\{x\}=
\pi_{{\sr 1},x}(T)$
and sometimes $T\{x\}\ne\pi_{1,x}(T).$ Returning to example 3.5, we get $T\{6\}=\pi_{\sr 1,6}(T)$ and
$T\{5\}\ne\pi_{\sr 1,5}(T),\ T\{3\}\ne \pi_{\sr 1,3}(T).$
\proclaim Proposition. For $T\in \bT_n^2$ let $T_2$ be the union of connected
subsequences $\{a_1,a_1+1,\ldots a_1+k_1\},\ldots, \{a_s,\ldots,a_s+k_s\}$
where $a_i>a_{i-1}+k_{i-1}+1$ for any $i\ :\ 1<i\leq s.$ Then
$$\Dscr(T)=\{T\<a_j+k_j>\ |\ 1\leq j\leq s\ {and}\ \pi_{1,a_j+k_j}(T)=T\{a_j+k_j\}\}.$$
\par
\Pf
By corollary 3.3 (ii) and proposition 3.15, one has $\Dscr(T)\st\{T\<s>\}_{s\in T_2}.$ Moreover, for any
$s\ : a_j\leq s< a_j+k_j$ one has $s\in \tau(T)$ and $s\not\in\tau(T\{s\}).$ Thus, $T\<s>\not>T.$
We obtain that $\Dscr(T)\st\{T\<a_j+k_j>\}_{j=1}^s.$ 
\par
Consider $T\pr=T\<a_j+k_j>.$ Since $T\pr_1=T_1+\{a_j+k_j\}$ and respectively
$\<T\pr_2>\st \<T_2>$ it is enough to show that the second condition of lemma 3.9 is satisfied, i.e.
for any $s\in T\pr_2$
one has $s^{\sr T\pr\{s\}}\ne a_j+k_j.$ Indeed, if $\pi_{1,a_j+k_j}(T)=T\{a_j+k_j\}$ one has that 
$s^{\sr T\pr\{s\}}>a_j+k_j$ for any $s>a_j+k_j$
 and $s^{\sr T\pr\{s\}}<a_j+k_j$ for any $s<a_j+k_j.$  
\par
On the other hand, if $T\{a_j+k_j\}\ne \pi_{1,a_j+k_j}(T)$ that means that for $a_{j+1}$ one has 
$a_{j+1}^{\sr T\{a_{j+1}\}}<a_j+k_j$
Thus, $a_{j+1}^{\sr T\pr\{a_{j+1}\}}=a_j+k_j.$ Hence, by lemma 3.9 $T\pr\not>T.$   
And this concludes the proof.
\QED
Again consider $T$ from example 3.5. One has
$$\Dscr
\left(\;
\vcenter{
\halign{& \hfill#\hfill
\tabskip4pt\cr
\multispan{5}{\hrulefill}\cr
\ss1
\vb & 1 &  &  3 &\ts\vb\cr\vs1
&&&&\cr
\ss1
\vb & 2 &  &  5 & \ts\vb\cr\vs1
&&&&\cr
\ss1
\vb & 4 &  &  6 & \ts\vb\cr\vs1
&&\multispan{3}{\hrulefill}\cr
\ss1
\vb & 7 &\ts\vb\cr\vs1
\multispan{3}{\hrulefill}\cr}}\right)=\left\{\;
\vcenter{
\halign{& \hfill#\hfill
\tabskip4pt\cr
\multispan{5}{\hrulefill}\cr
\ss1
\vb & 1 &  &  3 &\ts\vb\cr\vs1
&&&&\cr
\ss1
\vb & 2 &  &  5 & \ts\vb\cr\vs1
&&\multispan{3}{\hrulefill}\cr
\ss1
\vb & 4 & \ts\vb\cr\vs1
&&\cr
\ss1
\vb & 6 & \ts\vb\cr\vs1
&&\cr
\ss1
\vb & 7 &\ts\vb\cr\vs1
\multispan{3}{\hrulefill}\cr}}\right\}$$
\parno
{\bf 3.17}\ \ \ Finally, let us consider the case of tableaux with 
two rows. Let $(\bT_n^2)^\dagger$ denote the set of 
standard Young tableaux with two rows.  For any  
$S\in (\bT_n^2)^\dagger$ one has by 2.14 
$S=T^\dagger(w_{\sr S^\dagger})=T(\ov w_{\sr S^\dagger}).$
\par
By 1.7, for any $T,S\in (\bT_n^2)^\dagger$ on has $T\dos S$
(resp. $T\chos S$) iff $T^\dagger\dgs S^\dagger$ (resp. 
$T^\dagger\chgs S^\dagger$). 
\par
For any $S\in (\bT_n^2)^\dagger$ set $w_{\sr S}:=\ov w_{\sr S^\dagger}.$
By 2.2 (iii), for any $S,T\in (\bT_n^2)^\dagger$ one has $w_{\sr S}\dos w_{\sr T}$
iff $w_{\sr S^\dagger}\dgs w_{\sr T^\dagger},$ therefore, all the results for $\bT_n^2$
can be translated to $(\bT_n^2)^\dagger.$
\proclaim Theorem. \ \ Let $T,S\in (\bT_n^2)^\dagger.$ 
\item {(i)} One has $T\chos S$ iff $w_{\sr T}\dos w_{\sr S}.$
\item {(ii)} One has $T\dos S$ iff $w_{\sr T}\dos w_{\sr S}.$
\item {(iii)} $\Vscr_S\st\ov\Vscr_T$ iff $\gn\cap^{w_S}\gn \st \gn\cap^{w_T}\gn.$
\item{(iv)}   Orders $\do$ and $\cho$ coincide on $(\bT_n^2)^\dagger$.
\par

\bigskip
{\bf Acknowledgments.}\ \ \ The  problem  of combinatorial description
of inclusion of
orbital variety closures in terms of Young tableaux as well as the idea of 
the chain order and induced right
weak order on Young tableaux 
were suggested by A. Joseph. I would like to thank him for this and
for the fruitful discussions through various 
stages of this work.

I would  also like to express my gratitude to the referee.
His numerous remarks helped to improve the notation, alter some proofs
and bring this paper to a more orthodox, and, hopefully, readable form.

\bigskip
\parno
{\bf References}
\item{[D]}    {\twelverm M. Duflo,} Sur la classification des id\'eaux primitifs dans l'alg\`ebre
               enveloppante d'une alg\`ebre de Lie semi-simple, {\it Ann. of Math.} {\bf 105} (1977),
               107-130. 
\item{[F]}    {\twelverm W. Fulton,} ``Young tableaux'' LMSST 35, Cambridge University Press 1997.
\item{[H]}    {\twelverm W. Hesselink,} Singularities in the nilpotent 
                       scheme of a classical group, {\it Trans. Am. Math. 
                       Soc.} {\bf 222} (1976), 1-32.
\item{[J]}   {\twelverm A. Joseph,} Towards the Jantzen conjecture, Comp. Math. 40, 1980, pp 35-67.
\item{[JM]}   {\twelverm A. Joseph, A. Melnikov,} Quantization of hypersurface orbital varieties 
               in $\gs\gl_n$, The orbit method in geometry and physics. In honor of 
A.A. Kirillov, series ``Progress in Mathematics", 213, Birkhauser, 2003, 165-196.
\item{[Kn]}   {\twelverm D. E. Knuth,} ``The art of computer programming,'' 
                       Vol.3, Addison-Wesley (1969), 49-72.
\item{[M]}    {\twelverm A. Melnikov,} Robinson-Schensted 
procedure and combinatorial properties of geometric order in 
$\gs\gl(n),$
                   {\it C.R.A.S. I,} {\bf 315} (1992), 709-714.
\item{[M1]}  {\twelverm A. Melnikov,} On orbital variety 
closures in $\gs\gl_n.$ I. Induced Duflo order, J. of Algebra, 271, 2004, pp. 179-233.
\item{[Sa]}   {\twelverm B.E. Sagan,} The Symmetric Group, Graduate Texts in Mathematics 203, Springer 2000.
\item{[Sch]}  {\twelverm  M. P. Sch\"utzenberger,} La correspondance de Robinson, {\it LN in Math.} {\bf 597} (1976),
                59-113.
\item{[Sp1]} {\twelverm N.Spaltenstein,} The fixed point set of a unipotent 
                        transformation on the flag manifold, {\it Proc. Konin. 
                        Nederl. Akad.} {\bf 79} (1976), 452-456.
\item{[Sp2]}  {\twelverm N. Spaltenstein,} Classes unipotentes de 
                    sous-groupes de Borel, {\it LN in Math.} {\bf 964} (1982),
                     Springer-Verlag.
\item{[St1]}  {\twelverm R. Steinberg,} On the desingularization of the 
                        unipotent variety, {\it Invent. Math.} {\bf 36} 
                        (1976), 209-224.
\item{[St2]} {\twelverm R. Steinberg,} An occurrence of the 
                        Robinson-Schensted correspondence, {\it J. of Algebra}
                        {\bf 113} (1988), 523-528.

\bye